\documentclass[11pt,leqno]{amsart}
\usepackage[mathcal]{eucal}
\usepackage[utf8x]{inputenc}

\usepackage{amsmath,amssymb,amsthm}
\usepackage[dvipsnames]{xcolor}
\usepackage[colorlinks=true, linkcolor=MidnightBlue, citecolor=OliveGreen,pagebackref]{hyperref}
\usepackage[capitalise,nameinlink]{cleveref}
\usepackage{mathdots}

\newtheorem{theorem}{Theorem}[section]
\newtheorem{lemma}[theorem]{Lemma}
\newtheorem{proposition}[theorem]{Proposition}
\newtheorem{corollary}[theorem]{Corollary}

\theoremstyle{definition}
\newtheorem{definition}[theorem]{Definition}
\newtheorem{example}[theorem]{Example}

\numberwithin{equation}{section}

\DeclareMathOperator{\qund}{\quad\text{ and }\quad}

\DeclareMathOperator{\Aut}{Aut}
\DeclareMathOperator{\GL}{GL}
\DeclareMathOperator{\Mat}{Mat}
\DeclareMathOperator{\F}{\mathbb F}
\DeclareMathOperator{\Z}{\mathbb Z}
\DeclareMathOperator{\N}{\mathbb N}
\DeclareMathOperator{\Nor}{Nor}
\DeclareMathOperator{\Cen}{Cen}
\DeclareMathOperator{\Sym}{Sym}

\DeclareMathOperator{\id}{id}
\DeclareMathOperator{\Stab}{Stab}
\DeclareMathOperator{\stab}{stab}
\DeclareMathOperator{\reg}{reg}

\DeclareMathOperator{\rooted}{rt}
\DeclareMathOperator{\lay}{lay}
\DeclareMathOperator{\lab}{lab}

\newcommand{\dottimes}[1]{\times \overset{#1}{\dots} \times}
\newcommand*{\defeq}{\mathrel{\vcenter{\baselineskip0.5ex \lineskiplimit0pt\hbox{\scriptsize.}\hbox{\scriptsize.}}}=}
\newcommand{\GGS}{\textsc{GGS}}
\newcommand{\mGGS}{multi-\textup{GGS}}

\title[Automorphisms of multi-GGS groups]{The automorphism group of a multi-GGS group}

\author[J. M. Petschick]{Jan Moritz Petschick}
\address{Jan Moritz Petschick: Mathematisches
  Institut, Heinrich-Heine-Universit\"at, 40225 D\"usseldorf, Germany}
\email{jan.petschick@hhu.de}

\thanks{The research was funded by the Deutsche Forschungsgemeinschaft (DFG, German Research Foundation) — 380258175}


\begin{document}

\begin{abstract}
	A \mGGS-group is a group of automorphisms of a regular rooted tree, generalising the Gupta--Sidki~$p$-groups. We compute the automorphism groups of all non-constant \mGGS-groups.
\end{abstract}

\maketitle

\section{Introduction}

The family of Grigorchuk--Gupta--Sikdi-groups, hereafter abbreviated `\GGS-groups', is best known as a source of groups with exotic properties, e.g.\ just-infinite groups or infinite finitely generated periodic groups. It generalises earlier examples constructed by its three namesakes in the 80's, see~\cite{Gri80, GS83}. These groups are defined as groups of automorphisms of a $p$-regular rooted tree $X^*$, for an odd prime $p$. They are two-genereated, and one of the generators is defined according to a one-dimensional subspace $\mathbf{E} \subseteq \F_p^{p-1}$. Allowing $\mathbf{E}$ to be more than one-dimensional yields a natural generalisation, these `higher-dimensional' \GGS-groups are called \emph{\mGGS-groups} or \emph{multi-edge spinal groups}. We prefer the first term.

In many regards, \mGGS-groups do not differ overly much from their one-dimensional counterparts, e.g.\ they are periodic under similar conditions, see~\cite[Theorem~3.2]{AKT16}, they possess the congruence subgroup property, see~\cite{GU19}, and they allow similar branching structures. Their virtue, aside from extending the list of subgroups of $\Aut(X^*)$ with remarkable properties, lies therein that many conditions on \GGS-groups, when generalised to the higher-dimensional counterparts, reveal themselves as linear conditions. In this sense, \mGGS-groups are the more natural class.

We are concerned with the computation of the automorphism group of a given \mGGS-group. The automorphisms of groups acting on rooted trees have been investigated before, e.g.\ in~\cite{BS05,LN02}. In general, such groups are quite rigid objects, and their automorphisms are induced by homeomorphisms of the tree. Indeed, in many cases all automorphisms are actually induced by automorphisms of the tree, cf.~\cite{GW03,LN02}, and for some specific classes the automorphism groups can be uniformly computed, see~\cite{BS05}. More generally, the (abstract) commensurator of groups acting on rooted trees has been investigated, cf.~\cite{Rov02}; this is the group of `almost automorphisms', i.e.\ automorphisms between two finite-index subgroups.

However, there are only few explicit computations of the automorphism group of \GGS- and related groups. Sidki computed the automorphism group of the Gupta--Sidki $3$-group in \cite{Sid87}, and building on the approach for this group, the automorphism groups of the first Grigorchuk group \cite{GS04}, the Fabrikowski--Gupta and the constant \GGS-group on the ternary tree \cite{Sid87} have been computed. The first and the last two examples give a complete list of \GGS-groups acting on the ternary tree.

We now state our main result.

\begin{theorem}\label{AutmGGS:thm:main}
	Let $G$ be a non-constant \mGGS-group, and let $U$ be the maximal subgroup of $\F_p^\times$ such that $\mathbf E$ is invariant under the permutation action induced by $U$ by reordering the columns according to multiplication, and $W$ the maximal subgroup of $\F_p^\times$ of elements $\lambda$ such that $\mathbf{E} \subseteq \operatorname{Eig}_\lambda(u)$ for some $u \in U$. Then the following holds.
	\begin{enumerate}
		\item If $G$ is regular, then
		\[
			\Aut(G) \cong (G \rtimes \prod_\omega \mathrm{C}_{p})\rtimes (U \times W).
		\]
		\item If $G$ is symmetric, then
		\[
			\Aut(G) \cong (G \rtimes \mathrm{C}_p) \rtimes (U \times W).
		\]
	\end{enumerate}
\end{theorem}

The definitions of `regular' and `symmetric' can be found in \cref{sec:regular and symm}. The slightly obscure definitions of $U$ and $W$ are made more transparent in \cref{sec:coprime autom}.
We can immediately derive the following corollary.

\begin{corollary}
	Let $G$ be a non-constant \mGGS-group. Then the following statements hold.
	\begin{enumerate}
		\item The outer automorphism group of $G$ is finite if and only if $G$ is a symmetric \GGS-group.
		\item The outer automorphism group of $G$ is non-trivial.
		\item The automorphism group of $G$ contains elements of order coprime to $p$ if and only if~$\mathbf{E}$ is invariant under a permutation induced by multiplication in $\F_p$.
		\item The automorphism group of $G$ is a $p$-group if and only if $G$ is periodic and $\mathbf{E}$ is not invariant under any permutation induced by multiplication in $\F_p$.
	\end{enumerate}
\end{corollary}

We also explicitly compute the automorphism group for a selection of examples, e.g.\ all Gupta--Sidki $p$-groups, see \cref{sec:examples automggs}.

Our proof combines the methods developed by Sidki in~\cite{Sid87} (cf.\ \cite{BS05} for a sketch of the strategy used in Sidki's paper) with techniques used by the author to determine the isomorphism classes of \GGS-groups in~\cite{Pet19}. This, together with a theorem on the rigidity of branch groups of Grigorchuk and Wilson~\cite{GW03}, allows to reduce the complexity of the computations. On the other hand, the inclusion of the symmetric \GGS-groups complicates some of the arguments.

\section{Higher dimensional Grigorchuk--Gupta--Sidki-groups}

\subsection{Regular rooted trees and their automorphisms}
Let $p$ be an odd prime, and denote by $X$ the set~$\{0, \dots, p-1\}$. The Cayley graph~$X^*$ of the free monoid on $X$ is a $p$-regular rooted tree. We think of the vertices of~$X^*$ as words in~$X$. The root of the tree is the empty word $\varnothing$. We write $X^n$ for the set of all words of length~$n$, called the \emph{$n$-th layer of $X^*$}, and we identify $X$ and $X^1$.

Every (graph) automorphism $g \in \Aut(X^*)$ necessarily fixes the root, since it has a smaller valency than every other vertex. Consequently, every automorphism $g$ leaves all layers $X^n$ invariant. We write $\Stab(n)$ for the (setwise) stabilisier of $X^n$, and $\Stab_G(n)$ for its intersection with some subgroup $G \leq \Aut(X^*)$. We call a group $G \leq \Aut(X^*)$ \emph{spherically transitive} if it acts transitively on all layers $X^n$.

The group $\Aut(X^*)$ inherits the self-similar structure of $X^*$, and decomposes as a wreath product
\[
	\Aut(X^*) \cong \Aut(X^*) \wr_X \Sym(X).
\]
We deduce that $\Aut(X^*) \cong \Aut(X^*) \wr_{X^n} (\Sym(X) \wr \dottimes{p^n} \wr \Sym(X))$, for every $n \in \N$, where the finite iterated wreath product acts on $X^n$ as on the leaves of the the finite rooted $p$-regular tree with $n$ layers. The base group of the $n$\textsuperscript{th} such wreath product decomposition is equal to the $n$\textsuperscript{th} layer stabiliser. We denote the induced isomorphism $\Stab(n) \to \Aut(X^*) \dottimes{p^n} \Aut(X^*)$ by $\psi_n$. For $v \in X^n$, we denote the projection to the $v$\textsuperscript{th} component of the base group by $|_v \colon \Aut(X^*) \to \Aut(X^*)$, this so-called \emph{section map} is a group homomorphism on the pointwise stabiliser $\stab(v)$ of $v$. We call a subgroup $G \leq \Aut(X^*)$ \emph{self-similar} if $G|_v \subseteq G$ for all $v \in X^*$, and we call it \emph{fractal} if $\Stab_G(1)|_x \leq G$ for all $x \in X$.

The image of an element $g \in \Aut(X^*)$ in $\Sym(X)$ under the quotient by $\Stab(1)$ is denoted $g|^\varnothing$, and we write $g|^v = g|_v|^\varnothing$, for any $v \in X^*$, for the \emph{label of $g$ at $v$}. Any automorphism is uniquely determined by the collection of its labels.

We fix an embedding $\rooted \colon \Sym(X) \to \Aut(X^*)$ by $\rooted(\sigma)|^\varnothing = \sigma$ and $\rooted(\sigma)|^v = 1$ for all $v \in X^\ast \smallsetminus\{\varnothing\}$. We call the elements $\rooted{\Sym(X)}$ \emph{rooted automorphisms}.

Let $\Gamma \leq \Sym(X)$ be a permutation group. We define the \emph{$\Gamma$-labelled subgroup of $\Aut(X^*)$} by
\[
	\lab(\Gamma) = \{ g \in \Aut(X^*) \mid g|^v \in \Gamma \text{ for all } v \in X^* \}.
\]
It is a well-known fact that if $\Gamma$ is of order $p$, the subgroup $\lab(\Gamma)$ is a Sylow pro-$p$ subgroup of $\Aut(X^*)$.

Let $(x_i)_{i \in \N}$ be a sequence of elements $x_i \in X$. The words $\{x_0 \cdots x_k \mid k \in \N \}$ form a half-infinite \emph{ray} $R$ in $X^*$ (or, equivalenty, a point of the boundary). Write $\overline x$ for the ray associated to the constant sequence $(x)_{i \in \N}$. An \emph{$R$-directed~automorphism} is an automorphism $g$ fixing $R$ such that for all $v \in X^*$ either $v$ connected by an edge to an element of $R$, or
\[
	g|^{v} = \id.
\]
A spherically transitive group $G \leq \Aut(X^*)$ is \emph{regular branch over~$K$}, for a finite-index subgroup~$K \leq G$, if $K$ contains $\psi_1^{-1}(K \dottimes{p} K)$ as a subgroup of finite index.

\subsection{Multi-GGS-groups}\label{sec:regular and symm}
Fix the permutation $\sigma = (0 \, 1 \, \dots \, p-1)$. Write $a = \rooted(\sigma)$ and $A = \langle a \rangle$, as well as $\Sigma = \langle \sigma \rangle$. Let $\mathbf{E}$ be an $r$-dimensional subspace of $\F_p^{p-1}$, for $r > 0$. Choose an ordered basis $(\mathbf{b}_1, \dots, \mathbf{b}_r)$ of $\mathbf{E}$, and denote the standard basis of $\F_p^{r}$ by $(\mathbf{s}_1, \dots, \mathbf{s}_r)$. Let $E \in \Mat(r, p-1; \F_p)$ be the matrix with the basis elements as rows. The columns of $E$ are denoted $\mathbf{e}_i$ for $i \in \{1, \dots, p-1\}$; thinking of $\mathbf{E}$ as a subspace of $\{0\} \times \F_p^{p-1} \leq \F_p^p$, we will also write $\mathbf{e}_0$ for the zero column vector of length $r$. Define, for all $j \in \{1, \dots, r\}$, the $\overline{0}$-rooted automorphisms
\begin{align*}
	b^{\mathbf{s}_j} \defeq \psi_1^{-1}(b^{\mathbf{s}_j}, a^{\mathbf{s}_j \cdot \mathbf{e}_1}, \dots, a^{\mathbf{s}_j \cdot \mathbf{e}_{p-1}}).
\end{align*}
Since $a$ has order $p$, we may extend this definition to arbitrary vectors $\mathbf{n} \in \F_p^{r}$, such that
$\psi_1(b^{\mathbf{n}}) = (b^{\mathbf{n}}, a^{\mathbf{n}\cdot E})$, where for any $\mathbf{m} = (m_1, \dots, m_{p-1}) \in \F_p^{p-1}$ we set $a^{\mathbf{m}}$ to be the tuple $(a^{m_1}, \dots, a^{m_{p-1}})$ (and tuples are combined appropriately). The associated map $b^{\bullet}\colon \F_p^{r} \to \Aut(X^*)$ is an injective group homomorphism. We write $B$ for the image $b^{\F_p^{r}}$.

\begin{definition}
	The \emph{multi-\GGS-group associated to $\mathbf{E}$} is the group $G_\mathbf{E}$ of automorphisms generated by the set
	\[
		A \cup B.
	\]
\end{definition}
The subgroup $A$ (shared by all \mGGS-groups) is called the \emph{rooted group}, and the subgroup $B$ is called the \emph{directed group}. The generating set in the definition is clearly not minimal; a minimal generating set is given by $\{a\} \cup \{ b^{\mathbf{s}_j} \mid j \in \{1, \dots, r\} \}$.

If the dimension $r$ of $\mathbf{E}$ is $1$, one usually speaks of a \GGS-group rather than a \mGGS-group. In this case, abusing notation, we write $b$ for $b^{\mathbf{s}_1}$.

Depending on the space $\mathbf{E}$, we distinguish three classes of \mGGS-groups:
\begin{enumerate}
	\item If $\mathbf{E} = \{ (\lambda, \dots, \lambda) \in \F_p^{p-1} \mid \lambda \in \F_p \}$, we call $G_{\mathbf{E}}$ the \emph{constant \GGS-group}. This special case behaves very differently to all other \mGGS-groups; we will, for the most part, exclude it from our considerations.
	\item If $r=1$, the space $\mathbf{E}$ is contained in $\{ (\lambda_1, \dots, \lambda_{p-1}) \in \F_p^{p-1} \mid \lambda_i = \lambda_{p-i} \text{ for all } i \in \{ 1, \dots, p-1\} \}$, and $G_{\mathbf{E}}$ is not the constant \GGS-group, we call $G_{\mathbf{E}}$ a \emph{symmetric} \GGS-group.
	\item If $G_{\mathbf{E}}$ is neither constant nor symmetric, we call it a \emph{regular} \mGGS-group.
\end{enumerate}

We record some of the key properties of \mGGS-groups that have been established in the literature, cf.\ \cite[Proposition~3.3]{KT18} \& \cite[Lemma~2]{GU19} for statement (i), \cite[Proposition~4.3 and Proposition~3.2]{KT18} for statements (ii) and (iii), \cite[Lemma~3.5]{FZ13} for (iv), \cite[Theorem~C]{FGU17} for (v), and \cite[Proposition~3.1]{AKT16} for (vi).

\begin{theorem}\label{AutmGGS:thm:facts}
	Let $G = G_{\mathbf{E}}$ be a \mGGS-group. The the following statements hold.
	\begin{enumerate}
		\item If $G$ is regular, it is regular branch over the derived subgroup $G'$, and the equality $\psi_1(\Stab_G(1)') = G' \dottimes{p} G'$ holds.
		\item The abelisation of $G$ is an elementary abelian $p$-group of rank $r+1$.
		\item If $G$ is not constant, it is regular branch over $\gamma_3(G)$, such that $\psi_1(\gamma_3(\Stab_G(1))) = \gamma_3(G) \dottimes{p} \gamma_3(G)$.
		\item If $G$ is a symmetric \GGS-group, the intersection $\psi_1(G) \cap G' \dottimes{p} G'$ fulfils
		\[
			\psi_1(G) \cap (G' \dottimes{p} G') = \{ (g_0, \dots, g_{p-1}) \mid \textstyle\prod_{i = 0}^{p-1}g_i \in \gamma_3(G)\},
		\]
		and thus is of index $p$ in $G' \dottimes{p} G'$.
		\item Every \mGGS-group is self-similar and fractal.
	\end{enumerate}
\end{theorem}

It is fruitful to introduce the following overgroup to deal with the special case of symmetric \GGS-groups.

\begin{definition}
	Let $G$ be a \mGGS-group. Set
	\(
		\underline{c} = \psi_1^{-1}([b^{\mathbf{s}_1}, a], \id, \dots, \id).
	\)
	The \emph{regularisation $G_{\reg}$ of $G$} is the group
	\[
		G_{\reg} = \left\langle G \cup \{ \underline{c} \} \right\rangle.
	\]
\end{definition}

We record the following lemma on the regularisation of a \mGGS-group.

\begin{lemma}\label{AutmGGS:lem:facts implications for Greg}
	Let $G$ be a non-constant \mGGS-group. Then the following statements hold.
	\begin{enumerate}
		\item $G_{\reg} = G$ if and only if $G$ is regular.
		\item If $G$ is symmetric, then $|G_{\reg}: G| = p$.
		\item The derived subgroups of $G$ and $G_{\reg}$ are equal.
	\end{enumerate}
\end{lemma}

The first two statements are immediate consequences of \cref{AutmGGS:thm:facts}. Also the last statement follows, in view of
\[
	\psi_1([b^{\mathbf{s}_1}, \underline{c}]) = ([b^{\mathbf{s}_1}, [b, a]], \id, \dots, \id) \in \gamma_3(G) \dottimes{p} \gamma_3(G),
\]
and of
\[
	\psi_1([a, \underline{c}]) = ([b,a], \id, \dots, \id, [b,a]^{-1}) \in \{ (g_0, \dots, g_{p-1}) \mid \textstyle\prod_{i = 0}^{p-1}g_i \in \gamma_3(G)\},
\]
from \cref{AutmGGS:thm:facts}~(iii) and (iv).

\subsection{Constructions within \(\Aut(X^*)\)}

We introduce some notation. Let $g \in \Aut(X^\ast)$ and $n \in \N$. We define the \emph{$n$\textsuperscript{th}~diagonal of~$g$} as the element
\[
	\kappa_n(g) = \psi_n^{-1}(g, \dottimes{p^n}, g).
\]
Analogously, for any subset $G \subseteq \Aut(X^\ast)$ we define
\(
	\kappa_n(G) = \{ \kappa_n(g) \mid g \in G\}.
\)
Note that if $G$ is a group, the set $\kappa_n(G)$ is a group isomorphic to $G$.

\begin{definition}
	Let $S \subseteq \Aut(X^\ast)$ be a set of tree automorphisms. The \emph{diagonal~closure of~$S$} is the set
	\[
		\overline{S} = \left\{ \prod_{i = 0}^\infty \kappa_i(s_i) \;\middle|\; s_i \in S \text{ for } i \in \N \right\}.
	\]
\end{definition}
Since the $n$\textsuperscript{th} factor of the infinite product is contained in $\Stab(n)$, the product is defined as $\Aut(X^*)$ is closed in the (profinite) topology induced by the layer stabilisers. Note that the diagonal closure is in general not a subgroup, even if $S \leq \Aut(X^\ast)$ is one.

\begin{definition}
	Let
	$S = \rooted(\Sigma)$ be a group of rooted automorphisms of $\Aut(X^*)$. The group
	$$
		\kappa_{\infty}(S) = \langle \kappa_n(s) \mid n \in \N, s \in S \rangle
	$$
	is called the \emph{group of layerwise constant labels in $\Sigma$}.
\end{definition}
It is easy to see that $\kappa_{\infty}(S) \cong \prod_{\omega} S$, and $\overline{\kappa_{\infty}(S)} = \overline S$.

\subsection{Coordinates for \mGGS-groups} We first establish the following lemma, that allows us to uniquely describe elements of the first layer stabiliser in terms of `coordinates'. To be precise, we construct an isomorphism
\[
	\Stab_G(1) \cong (G' \dottimes{p} G') \rtimes B.
\]
This uses the fact that, also modulo $\psi_1^{-1}(G' \dottimes{p} G')$, the labels $g|^x$ at first layer vertices of an element $g \in \Stab_G(1)$ are completely determined by the image of $g$ in $G/G'$. Recall that $\mathbf{e}_i$ is the $i$\textsuperscript{th} column of $E$, and that $\mathbf{e}_0$ denotes the zero column vector of length $r$.

\begin{lemma}\label{AutmGGS:lem:global equations}
	Let $G$ be a non-constant \mGGS-group. Let $g_0, \dots, g_{p-1} \in G$ be a collection of elements of $G$. Then 
	\[
		\psi_1^{-1}(g_0, \dots, g_{p-1}) \in \Stab_{G_{\reg}}(1)
	\]
	if and only if there exist $\mathbf{n}_k \in \F_p^{r}$ and $y_k \in G'$ for $k \in \{0, \dots, p-1\}$ such that
	\[
		g_k = a^{s_k} b^\mathbf{n_k} y_k, \quad\text{where}\quad
		s_k = \sum_{i = 0}^{p-1} \mathbf{n}_i \cdot \mathbf{e}_{k-i}.
	\]
\end{lemma}

\begin{proof}
	We first prove that every element with its sections determined by a collection of vectors and elements of the commutator subgroup defines an element of the regularisation. Fix some $\mathbf{n}_k \in \F_p^r$ and $y_k \in G'$ for all $k \in X$. Then the element
	\[
		g = b^{\mathbf{n}_0}(b^{\mathbf{n}_1})^{a^{p-1}} \dots (b^{\mathbf{n}_{p-1}})^{a} \in \Stab_G(1)
	\]
	fulfils
	\begin{align*}
		g|_k &= a^{\mathbf{n}_{0} \cdot \mathbf{e}_{k}} a^{\mathbf{n}_{1} \cdot \mathbf{e}_{{k-1}}} \dots a^{\mathbf{n}_{k-1} \cdot \mathbf{e}_{{1}}} b^{\mathbf{n}_{k}} a^{\mathbf{n}_{k+1} \cdot \mathbf{e}_{{p-1}}} \dots a^{\mathbf{n}_{p-1} \cdot \mathbf{e}_{{k-1}}}\\
		&= a^{s_k} b^{\mathbf{n}_{k}} \tilde y_k
	\end{align*}
	for some $\tilde y_k \in G'$. We have
	\[
		\Stab_{G_{\reg}}(1) \geq \langle \Stab_G(1)' \cup \{ \underline{c} \} \rangle^G = \psi_1^{-1}(G' \times \dots \times G'),
	\]
	which follows directly from \cref{AutmGGS:thm:facts}~(i) for regular $G$. For symmetric $G$, by \cref{AutmGGS:thm:facts}~(iv), it is enough to show that $\langle \underline{c} \rangle^G = \psi_1^{-1}(G' \times \dots \times G')$. Clearly $[b,a]$ normally generates $G'$. The conjugates of $\underline{c}$ have only one non-trivial section, which is equal to $[b,a]$. The statement follows, since $G$ is fractal.
	
	Thus the element $y = \psi_1^{-1}(\tilde y_0^{-1}y_0, \dots, \tilde y_{p-1}^{-1}y_{p-1})$ is contained in $G_{\reg}$, so the element
	\[
		gy = \psi_1^{-1}(a^{s_0}b^{\mathbf{n}_{0}}y_0, \dots, a^{s_{p-1}}b^{\mathbf{n}_{p-1}}y_{p-1}) \in G_{\reg}
	\]
	has the prescribed sections.
	
	Now let $g \in \Stab_{G_{\reg}}(1)$. Up to $\psi_1^{-1}(G' \dottimes{p} G')$, i.e.\ up to the choice of $y_k \in G'$ for $k \in \{0, \dots, p-1\}$, we may calculate modulo the subgroup $L \defeq \langle \Stab_G(1)' \cup \{ \underline{c} \} \rangle^G \leq G_{\reg}$. Thus there are $\mathbf{n}_k \in \F_p^r$ for $k \in X$ such that
	\[
		g \equiv_L b^{\mathbf{n}_{0}}(b^{\mathbf{n}_{1}})^{a^{p-1}}\dots (b^{\mathbf{n}_{p-1}})^{a}.
	\]
	Taking sections as we did above shows that $g|_k \equiv_{G'} a^{s_k}b^{\mathbf{n}_{k}}$.
\end{proof}
Given $g \in \Stab_{G_{\reg}}(1)$, we call the vectors $\mathbf{n}_k$ introduced in \cref{AutmGGS:lem:global equations} the \emph{$B$-coordinates of~$g$}, and the collection of elements $y_k \in G'$ the \emph{$L$-coordinates~of~$g$}. The elements $s_k$ (since they are fixed by the $B$-coordinates) are called the \emph{forced~$A$-coordinates~of~$g$}.

\subsection{Strategy for the proof of \cref{AutmGGS:thm:main}}
By \cite[Theorem~1]{GW03} and \cite[Proposition~3.7]{KT18}, the automorphism group of $G$ coincides with the normaliser of $G$ in $\Aut(X^*)$. Hence we compute this normaliser $\Nor(G)$. In general, for any $H \leq \Aut(X^*)$, we denote by $\Nor(H)$ (without subscript) the normaliser of $H$ in $\Aut(X^*)$.

The normaliser of $\Sigma$ in $\Sym(X)$ has the form $\Nor_{\Sym(X)}(\Sigma) \cong \Sigma \rtimes \Delta$, where $\Delta \cong \F_p^\times$ with the multiplication action on $\Sigma \cong \F_p$. Heuristically, the automorphism group of a \mGGS-group $G$ allows for a similar decomposition. Since $G$ is contained in $\lab(\Sigma)$, its normaliser is contained in $\lab(\Nor_{\Sym(X)}(\Sigma))$, which decomposes as described above. The normaliser of $G$ within $\lab(\Sigma)$ is not identical to $G$, but turns out to be closely related. Apart from $G$ being symmetric or not, the structure of $\mathbf{E}$ only comes into play when considering the normaliser of $G$ in $\lab(\Delta)$.

We first consider the normaliser of $G$ in an appropriate closure within $\lab(\Sigma)$. Then we prove that the full normaliser splits as a semidirect product of the normaliser of $G$ within said closure, and normaliser of $G$ within an appropriate subgroup of $\lab(\Delta)$. At last, we compute the normaliser of $G$ within $\lab(\Delta)$, and combine our results.

\section{The normaliser in \(\overline{G_{\reg}}\)}
We begin our study of elements normalising $G$. Adating the strategy of Sidki in~\cite{Sid87}, we start not with the normaliser in the full automorphism group, but rather in the group $\overline{G_{\reg}} \geq G$. This a natural candidate, since it contains the normaliser of the rooted group $A$ (cf.\ \cref{AutmGGS:lem:norm a}) and the group $G$ itself.

\begin{lemma}\label{AutmGGS:lem:g infty diagonally closeds}
	Let $G$ be a non-constant \mGGS-group. Then
	\[
		\kappa_1(G_{\reg}) \leq \kappa_{\infty}(A) \cdot G_{\reg}.
	\]
\end{lemma}

\begin{proof}
	We check that the generators of $\kappa_1(G_{\reg})$ are contained in the group on the right hand side. Clearly $\kappa_1(a) \in \kappa_{\infty}(A)$.
	
	To see that $\kappa_1(b^{\mathbf{s}_j})$ is contained in $\kappa_{\infty}(A) \cdot G_{\reg}$ for a given $j \in \{ 1, \dots, r\}$, we use \cref{AutmGGS:lem:global equations}.
	We have no choice for the set of $B$-coordinates; since $\kappa_1(b^{\mathbf{s}_j})|_x = b^{\mathbf{s}_j}$ for all $x \in X$ they are all equal to $\mathbf{s}_j$. Thus we compute the forced $A$-coordinates
	\[
		s_k = \sum_{i = 0}^{p-1} \mathbf{n}_i \cdot \mathbf{e}_{{k-i}} = \sum_{i = 0}^{p-1}\mathbf{s}_j \cdot \mathbf{e}_{{k-i}} = \sum_{i = 0}^{p-1} e_{j, k-i},
	\]
	where $e_{j, k-i}$ is the respective entry of $E$. Consequently, all forced $A$-coordinates are equal to some fixed $s \in \F_p$ and independent of $k$. Hence
	\[
		\kappa_1(a^s b^{\mathbf{s}_j}) \in G_{\reg}.
	\]
	Since we have already established that $\kappa_1(a) \in \kappa_{\infty}(A)$, this implies that $\kappa_1(b^{\mathbf{s}_j})$ is contained in $\kappa_{\infty}(A) \cdot G_{\reg}$ for all $j \in \{1, \dots, r\}$.
	
	Finally, in the case that $G$ is symmetric, we have $[\kappa_1(a), b] = \underline{c}$, and hence
	\[
		\psi_1([\kappa_2(a), \kappa_1(b)]) = \kappa_1([\kappa_1(a), b]) = \kappa_1(\underline{c}).\qedhere
	\]
\end{proof}

The problem to determine the normaliser is easily solved for the rooted group $A$. To determine the normaliser of $B$ is significantly harder.

\begin{lemma}\label{AutmGGS:lem:norm a}
	The centraliser and the normaliser of the rooted group $A$ are given by
	\begin{align*}
		\Cen(A) &= \kappa_1(\Aut(X^*)) \rtimes A, \quad \text{and}\\
		\Nor(A) &= \kappa_1(\Aut(X^*)) \rtimes \rooted(\Nor_{\Sym(X)}(\Sigma)).
	\end{align*}
\end{lemma}

\begin{proof}
	Given $x \in X$, we have \(a^g|_x = \id\) if and only if $g|_{x} = g|_{x+1}$, hence we have $a^g|_x = a^j|_x$ for some $j \in \{1, \dots, p-1\}$ if and only if $g|_0 = g|_x$ for all $x \in X$. The image of $a$ under conjugation with $g$ now only depends on $g|^\varnothing$, hence we only need to observe $\Cen_{\Sym(X)}(\Sigma) = \Sigma$.
\end{proof}

\begin{lemma}\label{AutmGGS:lem:centraliser and normaliser of b restricts to itself in section at 0}
	Let $G$ be a non-constant \mGGS-group. Then $\Nor(B) \leq \stab(0)$, the point stabiliser of the vertex $0 \in X$, and
	\[
		\Nor(B)|_0 \leq \Nor(B) \qund \Cen(B)|_0 \leq \Cen(B).
	\]
\end{lemma}

\begin{proof}
	Let $g \in \Nor(B)$. Then there is some $\mathbf{n} \in \F_p^r\smallsetminus\{0\}$ such that $(b^{\mathbf{s}_1})^g = b^{\mathbf{n}}$. If $0^{g^{-1}} = x \neq 0$, we see that
	\[
		b^{\mathbf{n}} = b^{\mathbf{n}}|_0 = (b^{\mathbf{s}_1})^g|_0 = (a^{\mathbf{b}_{i, x}})^{g|_0}.
	\]
	But a rooted automorphism cannot be conjugate to a directed automorphism. Thus $g \in \stab(0)$. Similarly, we find
	\[
		b^{\mathbf{n}} = (b^{\mathbf{s}_1})^g|_0 = (b^{\mathbf{s}_1}|_0)^{g|_0} = (b^{\mathbf{s}_1})^{g|_0}.
	\]
	This shows (allowing $\mathbf{n} = \mathbf{s}_1$) both other statements.
\end{proof}

For the next lemma we introduce the (word) length function $\ell \colon G \to \N$, with respect to the generating set $A \cup B$, i.e.\ the mapping
\[
	\ell(g) = \min\{n \in \N \mid g \text{ can be written as a product of length $n$ in } A \cup B\}.
\]
It is well-known that this length function is \emph{contracting}, i.e. that $\ell(g|_x) \leq g$ for $x \in X$. We need some finer analysis to establish a strict inequality in certain cases. Note that the strictness of the inequality above, for a more general class of self-similar groups $G$, is related to $G$ being a periodic group.

\begin{lemma}\label{AutmGGS:lem:length reduction alternating}
	Let $G$ be a non-constant \mGGS-group, and let $g \in G$ be an element with $\ell(g) > 1$. Then there is some $i \in X\smallsetminus\{0\}$ such that $\ell(g|_0g|_i^{-1}) < \ell(g)$.
\end{lemma}

\begin{proof}
	Write $g = a^{i_0} b^{\mathbf{n}_{0}} \dots a^{i_{n-1}} b^{\mathbf{n}_{m-1}} a^{i_m}$, where $m \in \N$, $\mathbf{n}_{k} \in \F_p^\times$, $i_{k} \in \Z\smallsetminus\{0\}$ for $k \in  \{0, \dots, m-1\}$, and $i_{n} \in \Z$. Passing to a conjugate if necessary, every $g \in G$ can be written in this way. Taking sections, we see that
	\[
		g|_k = b^{\mathbf{n}_{0}}|_{k-i_0} (b^{\mathbf{n}_{1}})|_{k-i_0-i_1} \dots (b^{\mathbf{n}_{m-1}})|_{k - \sum_{t = 0}^{n-1}i_t}
	\]
	for any $k \in X$, hence $\ell(g|_k) \leq m$. Since every $B$-letter contributes at most one $B$-letter to one of the sections, we have $\sum_{k=0}^{p-1} \ell(g|_k) \leq \ell(g) + p - 1$. Assume that $\ell(g|_0g|_1^{-1}) \geq \ell(g)$. Then $\ell(g|_0) + \ell(g|_1) \geq \ell(g)$, hence $\ell(g|_0) = \ell(g|_1) = m$. This can only be the case if every $B$-letter contributes its only section that is contained in $B$ either to $g|_0$ and $g|_1$, i.e.\ $g|_k \in A$ for all other $k \in X$. Thus, if $m > 2$, we have
	\[
		\ell(g|_0g|_k^{-1}) \leq \ell(g|_0) + \ell(g|_k) \leq m + 1 \leq 2m - 1 < \ell(g),
	\]
	If $m = 2$, we see that $g = a^{i_0}b^{\mathbf{n}_{0}}$, implying $g|_0 \in A$. Since there is at least a second section contained in $A$, the result follows.
\end{proof}

\begin{lemma}\label{AutmGGS:lem:Cen a into G1 basic}
	Let $G$ be a non-constant \mGGS-group. Then
	\[
		\Nor_{\Cen(A)}(G) \cap (\Cen(B)\cdot G) \subseteq \kappa_{\infty}(A) \cdot G_{\reg}.
	\]
\end{lemma}

\begin{proof}
	Let $g \in \Nor_{\Cen(A)}(G) \cap (\Cen(B)\cdot G)$ and $h \in G$ be an element of minimal length such that we may write $g = g' h$ for some $g' \in \Cen(B)$. The proof uses induction on the length of $h$.
	
	First assume that $h$ has length one, i.e.\ $h \in A \cup B$. If $h$ is in $B$, we find that $g \in \Cen(B)$. Thus $h$ centralises $G$, but it is well-known that the centraliser of a branch group in $\Aut(X^*)$ is trivial; hence $g = \id$. If $h$ is a power of $a$, the same holds for $gh^{-1}$, hence $g \in A \leq G_{\reg}$.
	
	Now we assume that $\ell(h) > 1$. By \cref{AutmGGS:lem:norm a} we may write $g = \kappa_1(g|_0)a^k$ for some $k \in \Z$, yielding for any $\mathbf n \in \F_p^r$
	\begin{align*}
		((b^{\mathbf{n}})^{g|_0}, (a^{\mathbf n \cdot \mathbf{e}_{{1}}})^{g|_0}, \dots, (a^{\mathbf n \cdot \mathbf{e}_{{p-1}}})^{g|_0})^{a^k} &= \psi_1((b^{\mathbf{n}})^g) = \psi_1((b^{\mathbf{n}})^h)\\
		&= ((b^{\mathbf n})^{h|_0}, (a^{\mathbf n \cdot \mathbf{e}_{{1}}})^{h|_1}, \dots, (a^{\mathbf n \cdot \mathbf{e}_{{p-1}}})^{h|_{p-1}})^{h|^\varnothing}.
	\end{align*}
	Since $a$ and $b^{\mathbf n}$ are not conjugate in $\Aut(X^\ast)$, this shows that $g|^\varnothing = h|^\varnothing$ and $a^{g|_0} = a^{h|_i}$ for all $i \in X\smallsetminus\{0\}$. Thus $g|_0h|_i^{-1}$ centralises $A$, and by \cref{AutmGGS:lem:centraliser and normaliser of b restricts to itself in section at 0} we find $(b^{\mathbf n})^{g|_0h|_i^{-1}} = (b^{\mathbf n})^{h|_0h|_i^{-1}} \in B^G$, hence $g|_0h|_i^{-1}$ normalises $G$. By \cref{AutmGGS:lem:length reduction alternating}, there is some $i \in X\smallsetminus\{0\}$ such that $\ell(h|_0h|_i^{-1}) < \ell(h)$, so by induction we see that $g|_0h|_i^{-1} \in \kappa_{\infty}(A) \cdot G_{\reg}$. Since $h|_i \in G$, we have $g|_0 \in \kappa_{\infty}(A) \cdot G_{\reg}$, and
	$$
		g = \kappa_1(g|_0) a^{k} = a^k \kappa_1(g|_0) \in A \cdot \kappa_1(\kappa_{\infty}(A) \cdot G_{\reg}) = \kappa_{\infty}(A) \cdot \kappa_1(G_{\reg}).
	$$
	Now \cref{AutmGGS:lem:g infty diagonally closeds} yields $g \in \kappa_{\infty}(A) \cdot G_{\reg}$.
\end{proof}

With a little care, we can use the same idea to extend the result to $G_{\reg}$.

\begin{lemma}\label{AutmGGS:lem:Cen a into G1}
	Let $G$ be a non-constant \mGGS-group. Then
	$$\Nor_{\Cen(a)}(G) \cap (\Cen(b) \cdot G_{\reg}) \leq \kappa_{\infty}(A) \cdot G_{\reg}.$$
\end{lemma}

\begin{proof}
	In view of the previous lemma, we may restrict to symmetric $G$. Let $g \in \Nor_{\Cen(A)}(G) \cap (\Cen(B) \cdot G_{\reg})$ and choose $g'\in \Cen(B)$, $h \in G$ and $j \in \Z$, such that $g = g'\underline{c}^jh$. 
	Write $g = \kappa_1(g|_0)a^k$ for some $k \in \Z$, and calculate,
	\begin{align*}
		(b^{g|_0}, (a^{\mathbf{e}_{{1}}})^{g|_0}, \dots, (a^{\mathbf{e}_{{p-1}}})^{g|_0})^{a^k} &= \psi_1(b^g) = \psi_1(b^{\underline{c}^jh}) \\&= (b^{\underline{c}^j|_0h|_0}, (a^{\mathbf{e}_{{1}}})^{h|_1}, \dots, (a^{\mathbf{e}_{{p-1}}})^{h|_{p-1}})^{h|^\varnothing}.
	\end{align*}
	As we did in the proof of \cref{AutmGGS:lem:Cen a into G1 basic}, we may conclude that $h|^\varnothing = a^k$. Consequently, for all $i \in X\smallsetminus\{0\}$, the element $g|_0h|_j^{-1}$ centralises $a$.
	By \cref{AutmGGS:lem:centraliser and normaliser of b restricts to itself in section at 0}, the element $g|_0$ is in $\Cen(B) \cdot \underline{c}^j|_0 h|_0$. Since $\underline{c}^j|_0 = [b,a]^j \in G$, this implies that $g|_0h|_i^{-1}$, for all $i \in X\smallsetminus\{0\}$, is an element in $\Nor_{\Cen(A)}(G) \cap (\Cen(B) \cdot G)$. By \cref{AutmGGS:lem:Cen a into G1 basic}, the element $g|_0h|_i^{-1}$, and consequently also $g|_0$ is contained in $\kappa_{\infty}(A) \cdot G_{\reg}$. Finally, by \cref{AutmGGS:lem:g infty diagonally closeds}, we find $g = \kappa_1(g|_0)g|^\varnothing \in \kappa_{\infty}(A) \cdot G_{\reg}$.
\end{proof}

\begin{lemma}\label{AutmGGS:lem:norm = cent} Let $G$ be a non-constant \mGGS-group. Then
	\begin{align*}
		\Nor_{\operatorname{lab}(\Sigma)}(A) \leq \Cen(A) \qund
		\Nor_{\operatorname{lab}(\Sigma)}(B) \leq \Cen(B).
	\end{align*}
\end{lemma}

\begin{proof}
	We use the description of $\Nor(A)$ given in \cref{AutmGGS:lem:norm a}. Let $g \in \Nor_{\lab(\Sigma)}(A)$. For all $h \in \lab(\Sigma)$, we have $h|^\varnothing \in \langle \sigma \rangle$. Thus we see that $a^{\rooted\kappa_1(g|_0)g|^\varnothing} = a^{\rooted g|^\varnothing} = a$.
	
	Now let $g \in \Nor_{\lab(\Sigma)}(B)$, let $\mathbf n \in \F_p^r$ be arbitrary and let $\mathbf{m} \in \F_p^r$ be such that $(b^{\mathbf{n}})^g = b^{\mathbf{m}}$. Then
	\[
		(b^{\mathbf{m}}, a^{\mathbf m \cdot \mathbf{e}_{{1}}}, \dots, a^{\mathbf m \cdot \mathbf{e}_{{p-1}}}) = b^{\mathbf m} = (b^{\mathbf n})^g = ((b^{\mathbf n})^{g|_0}, (a^{\mathbf n \cdot \mathbf{e}_{{1}}})^{g|_1}, \dots, (a^{\mathbf n \cdot \mathbf{e}_{{p-1}}})^{g|_{p-1}})^{g|^\varnothing}.
	\]
	The label $g|^\varnothing$ is a power of $a|^\varnothing$. Since $b^{\mathbf{n}}$ and $a$ are not conjugate in $\Aut(X^\ast)$, the element $g|^\varnothing$ must stabilise the vertex $0$, thus it is trivial. Varying $\mathbf{n}$, we see that $g|_i$ normalises $A$ for all $i \in X\smallsetminus\{0\}$ for which $\mathbf{e}_{{i}} \neq \mathbf{0}$.
	Now since $\lab(\Sigma)$ is self-similar, this implies $g|_i \in \Cen(A)$ by the first part of this lemma, hence $a^{\mathbf n \cdot \mathbf{e}_{{i}}} = a^{\mathbf m \cdot \mathbf{e}_{{i}}}$ for all $i \in X\smallsetminus\{0\}$. Thus $b^{\mathbf m} = b^{\mathbf n}$.
	Since $g|_0 \in \Nor_{\lab(\Sigma)}(B)$ by \cref{AutmGGS:lem:centraliser and normaliser of b restricts to itself in section at 0}, we can argue in the same way for $g|_0$, hence $g \in \Cen(B)$.
\end{proof}

\begin{lemma}\label{AutmGGS:lem:finite words in g1}
	Let $G$ be a non-constant \mGGS-group. Then
	$$\Nor_{\overline{G_{\reg}}}(G) \subseteq \kappa_{\infty}(A)\cdot G_{\reg}.$$
\end{lemma}

\begin{proof}
	Let $g \in \Nor_{\overline{G_{\reg}}}(G)$. There is a sequence $(g_i)_{i\in \N}$ with $g_i \in G_{\reg}$ such that
	\[
		g = \prod_{i = 0}^{\infty} \kappa_i(g_i).
	\]
	Write $h_n$ for the partial product $\prod_{i = 0}^{n} \kappa_i(g_i)$. By \cref{AutmGGS:lem:g infty diagonally closeds}, we find $h_n \in \kappa_{\infty}(A) \cdot G_{\reg}$ for all $n \in \N$. We may write
	\begin{align*}
		g|_{0^n} = h_n|_{0^n} \left(\prod_{i = n+1}^{\infty} \kappa_i(g_i)\right)|_{0^n.h_n}
		= h_n|_{0^n} \prod_{i = 1}^{\infty} \kappa_i(g_{i+n}).
	\end{align*}
	In view of \cref{AutmGGS:lem:norm a}, we conclude that $h_n|_{0^n}^{-1}g|_{0^n} \in \Cen(A)$. There is nothing special about $0^n$; indeed, we see that $h_n|_v^{-1}g|_v = h_n|_{0^n}^{-1}g|_{0^n}$ for all $v \in X^n$. 
	By \cite[Lemma 3.4]{Pet19}, there exists an integer $n \in \N$ such that $g|_{0^n} \in \Nor(B)$. By \cref{AutmGGS:lem:norm = cent} $g|_{0^n} \in \Cen(B)$. Consequently
	\[
		B^{g|_{0^n}^{-1}h_n|_{0^n}} = B^{h_n|_{0^n}}.
	\]
	Since $h_n|_{v} \in G_{\reg}$ for all $v \in X^n$, we may use \cref{AutmGGS:lem:Cen a into G1} and obtain $g|_{0^n}^{-1}h_n|_{0^n} \in \kappa_{\infty}(A) \cdot G_{\reg}$. Using \cref{AutmGGS:lem:g infty diagonally closeds} again, we find $\kappa_n(h_n|_{0^n}^{-1}g|_{0^n}) \in \kappa_{\infty}(A) \cdot G_{\reg}$, and moreover
	\[
		g = h_n \psi_n^{-1}(h_n|_{0^n}^{-1}g|_{0^n}, \dots, h_n|_{(p-1)^n}^{-1}g|_{(p-1)^n}) = h_n \kappa_n(h_n|_{0^n}^{-1}g|_{0^n}) \in \kappa_{\infty}(A) \cdot G_{\reg}.\qedhere
	\]
\end{proof}

\begin{lemma}\label{AutmGGS:lem:glay and its normalising part}
	Let $G$ be a non-constant \mGGS-group. Write $G_{\lay}$ for the product set $\kappa_{\infty}(A) \cdot G_{\reg}$.
	\begin{enumerate}
		\item If $G$ is regular, then we have $\kappa_{\infty}(A) \leq \Nor_{\Aut(X^*)}(G)$, hence $G_{\lay}$ acquires the structure of a semidirect product.
		\item If $G$ is symmetric, we find $\Nor_{\kappa_{\infty}(A)}(G) = A$.
		\item Write $G_{\lay}$ for product set in \textup{(i)}. Then
		\[
			\Nor_{G_{\lay}}(G) = \begin{cases}
				G_{\lay} &\text{ if }G\text{ is regular},\\
				G_{\reg} &\text{ if }G\text{ is symmetric.}\\
			\end{cases}
		\]
	\end{enumerate}
\end{lemma}

\begin{proof}
	Let $n \in \N$. Clearly $a^{\kappa_n(a)} = a$, and for all $j \in  \{1, \dots, r \}$
	\begin{align*}
		[b^{\mathbf{s}_j},\kappa_n(a)] &= \psi_1^{-1}([b^{\mathbf{s}_j}, \kappa_{n-1}(a)], [a^{\mathbf{b}_{j,1}}, \kappa_{n-1}(a)], \dots, [a^{\mathbf{b}_{j,p-1}}, \kappa_{n-1}(a)])\\
		&= \psi_n^{-1}([b^{\mathbf{s}_j},a], \id, \dots, \id) \in \psi_n^{-1}(G' \times \dots \times G') \leq G.
	\end{align*}
	This shows (i), and it also shows that $\kappa_n(a)$ does not normalise a symmetric \GGS-group $G$ for $n > 0$, since $\psi_1^{-1}([b, a], \id, \dots, \id) \notin G$. Thus (ii) is proven.
	
	Statement (iii) is a consequence of (i) in case $G$ is regular, and an immediate consequence of \cref{AutmGGS:lem:facts implications for Greg}~(iii) in case $G$ is symmetric.
\end{proof}

\begin{proposition}\label{AutmGGS:prop:norm in g1}
	Let $G$ be a non-constant \mGGS-group. Then
	\[
		\Nor_{\overline{G_{\reg}}}(G) = \begin{cases}
			G \rtimes \kappa_{\infty}(A), &\text{ if $G$ is regular},\\
			G_{\reg}, &\text{ if $G$ is symmetric.}
		\end{cases}
	\]
\end{proposition}

\begin{proof}
	Assume that $G$ is regular. By \cref{AutmGGS:lem:glay and its normalising part}, the set $\kappa_{\infty}(A) \cdot G_{\reg}$ is a group. In view of \cref{AutmGGS:lem:finite words in g1} and $G = G_{\reg}$, this proves the first case. If $G$ is symmetric, the result follows from \cref{AutmGGS:lem:finite words in g1} and \cref{AutmGGS:lem:glay and its normalising part}.
\end{proof}

\section{The normaliser as a product}

We now prove that the normaliser of $G$ in $\Aut(X^*)$ decomposes as a semi-direct product. To begin with, we prove the following generalisation of \cite[2.2.5(i)]{Sid87}, which is an interesting proposition in its own right.

\begin{proposition}\label{AutmGGS:prop:order p elements}
	Let $G$ be a non-constant \mGGS-group. Every element of $G$ that has order $p$ is either contained in $\Stab_G(1)$ or is conjugate to a power of $a$ in $G_{\reg}$.
\end{proposition}

\begin{proof}
	We have to prove that, given $g \in \Stab_G(1)$ and $i \in \Z$, every element $a^ig$ of order $p$ may be written $(a^h)^i$ for some $h \in G_{\reg}$. Passing to an appropriate power of $a$, we may assume that $i = 1$.
	From $(ag)^p = 1$ we derive the equations
	\begin{align*}
		\id &= (ag)^p|_0 = g|_{0} \dots g|_{p-1}, \text{ resp.}\\
		g|_{p-1} &= g|_{p-2}^{-1} \dots g|_0^{-1}.
	\end{align*}	
	Since $g \in \Stab_G(1)$,  by \cref{AutmGGS:lem:global equations} there exists a set of $B$-coordinates $\mathbf{n}_{k} \in \F_p^r$ and a set of $L$-coordinates $y_k \in G'$ uniquely describing $g$. Reformulated in these $B$-coordinates, the condition above reads
	\[
		\sum_{i = 0}^{p-2} \mathbf{n}_i = -\mathbf{n}_{p-1}.
	\]
	Given some integer $s \in \Z$, we define an element
	\begin{align*}
		h_s = \psi_1^{-1}(a^s, a^sg|_0, a^sg|_0g|_1, \dots, a^sg|_0 \dots g|_{p-2}) \in \psi_1^{-1}(G \dottimes{p} G).
	\end{align*}
	Since
	\begin{align*}
		a^{h_s}|_k &= h_s|_{k}^{-1}h_s|_{k+1} = (a^s\prod_{i = 0}^{k-1}g|_i)^{-1}a^s\prod_{i = 0}^{k}g|_i \\&= \begin{cases}
			g|_k, &\text{ if }k \neq p-1,\\
			(\prod_{i = 0}^{p-2}g|_i)^{-1} = g|_{p-1}, &\text{ if }k = p-1,\\
		\end{cases}
	\end{align*}
	the conjugate $a^{h_s}$ is equal to $ag$. It remains to prove that $h_s \in G_{\reg}$ for some $s \in \Z$. If it is contained in $G_{\reg}$, the element $h_s$ has the $B$-coordinates ${\mathbf{h}_k} = \sum_{i = 0}^{k-1} \mathbf{n}_i$ (and some commutators $z_k$ that we shall not need to specify). We have to prove that the corresponding forced $A$-coordinates $\tilde s_k = \sum_{i = 0}^{p-1} \mathbf{h}_i \cdot \mathbf{e}_{{k-i}}$ are equal to the actual $a$-exponents of the corresponding sections of $h$. Since it is enought to show that $h_s \in G_{\reg}$ for one $s$, we fix $s = \tilde s_0$, so that the proposed equality holds in the first component by definition.
	
	A quick calculation shows that, for all $k \in X \smallsetminus \{0\}$,
	\begin{align*}
		\tilde s_{k} - \tilde s_{k-1} &= \sum_{i = 0}^{p-1} \mathbf{h}_i \cdot \mathbf{e}_{{k-i}} - \sum_{i = 0}^{p-1} \mathbf{h}_i \cdot \mathbf{e}_{{k-1-i}}\\
		&= \sum_{i = 0}^{p-1} (\mathbf{h}_{i} - \mathbf{h}_{i-1} ) \cdot \mathbf{e}_{{k - i}}
		= \sum_{i = 0}^{p-1} \mathbf{n}_{i-1} \cdot \mathbf{e}_{{k-i}} = s_{k-1},
	\end{align*}
	and consequently the $a$-exponent of $h|_k$ is equal to
	\[
		s + \sum_{i = 0}^{k-1} s_i = \tilde s_0 + \sum_{i = 0}^{k-1} s_i = \tilde s_k,
	\]
	for all $k \in X$. But the values $s_i$ are the forced $A$-coordinates of $g$, hence, comparing with the definition of $h_s$, we see that the forced $A$-coordinates of $\mathbf{h}_k$, for $k \in X$, and the actual $a$-exponents of $h|_k$ coincide. Hence $h \in G_{\reg}$.
\end{proof}

Notice that for a symmetric \GGS-group, we do have to pass to $G_{\reg}$ to make this statement true: take the element $d = ([b, a], [b,a]^{-1}, \id, \dots, \id) \in G$. Clearly $a^{[b, a]} = d$, but assume for contradiction that there is another element $h \in \Stab_G(1)$ such that $a^{h^{-1}} = d$. Then $\underline{c}h$ centralises $a$, hence
\[
	[b,a]h|_0 = (\underline{c}h)|_0 = (\underline{c}h)|_i = h|_i,
\]
for all $i \in \F_p^\times$. Counting the powers of $[b,a]$ in $h|_k$ as in \cref{AutmGGS:thm:facts}~(iv), we see that if $h|_1 \equiv_{\gamma_3(G)} a^s(b)^n([b,a])^v$ the sum of the $[b,a]$-exponents over all sections mod $\gamma_3(G)$ equals $v-1 + (p-1)v \equiv_p p-1$, contradicting \cref{AutmGGS:thm:facts}~(iv). Thus there is no such $h \in \Stab_G(1)$.

Recall that the group $\Delta$ is $\Nor_{\Sym(X)}(\Sigma) \cap \stab_{\Sym(X)}(0)$. Set $D = \rooted(\Delta)$, i.e.\ the group of rooted automorphisms normalising but not centralising $a$.

\begin{lemma}\label{AutmGGS:lem:norm in ginfty kappa closure}
	Let $G$ be a non-constant \mGGS-group. Then
	\[
		\Nor(G) \subseteq \overline{G_{\reg}} \cdot \overline{D}.
	\]
\end{lemma}

\begin{proof}
	Let $g_0 \in \Nor(G)$. Let $k \in \Z$ be such that $(a^{g_0})^k|^\varnothing = a$. By \cref{AutmGGS:prop:order p elements} there exists an element $h_0 \in G_{\reg}$ such that
	\(
		(a^{g_0})^k = a^{h_0}.
	\)
	Consequently $h_0^{-1}g_0 \in \Nor(A)$. Using \cref{AutmGGS:lem:norm a} and the fact that $\Nor(G)$ is self-similar, cf.~\cite[Lemma~3.3]{Pet19}, we may write
	\[
		h_0^{-1}g_0 = \kappa_1((h_0^{-1}g_0)|_0)\rooted((h_0^{-1}g_0)|^\varnothing)
	\]
	for $h_0^{-1}g_0|_0 = g_1 \in \Nor(G)$. Since $g_0h_0^{-1}|^\varnothing$ normalises $\sigma$, we may write $g_0h_0^{-1}|^\varnothing = a^{k_0}d_0$ for some $d_0 \in D$ and $k_0 \in \Z$, so that we obtain the equation
	\[
		g_0 = h_0\kappa_1(g_1) a^{k_0} d_0 = h_0 a^{k_0} \kappa_1(g_1) d_0,
	\]
	using the fact that $\kappa_1(\Aut(X^*))$ normalises $a$ for the second equality. Repeating the procedure for $g_1$, we obtain $g_2 \in \Nor(G), h_1 \in G_{\reg}, d_1 \in D$ and $k_1 \in \Z$ such that
	\begin{align*}
		g_0 &= h_0 a^{k_0} \kappa_1(h_1 a^{k_1} \kappa_1(g_2) d_1) d_0\\
		&= h_0 a^{k_0} \kappa_1(h_1 a^{k_1}) \kappa_2(g_2) \kappa_1(d_1) d_0\\
		&= h_0 a^{k_0} \kappa_1(h_1 a^{k_1}) \kappa_2(g_2) d_0 \kappa_1(d_1).
	\end{align*}
	In the last step we have used the fact that $\overline{D}$ is abelian.
	Going on, we obtain a sequence of products
	\[
		t_n = \prod_{i = 0}^n \kappa_{i}(h_ia^{k_i})
		\prod_{i = 0}^n \kappa_{n-i}(d_i) = \prod_{i = 0}^n \kappa_{i}(h_ia^{k_i})
		\prod_{i = 0}^n \kappa_{i}(d_i),
	\]
	such that $t_n \equiv_{\Stab(n+1)} g_0$, i.e.\ that are converging to $g_0$ in the topology induced by the layer stabilisers. Since both $\overline{D}$ and $\overline{G_{\reg}}$ are closed sets, the corresponding limites are well-defined. We obtain
	\[
		g_0 = \prod_{i = 0}^\infty \kappa_{i}(h_ia^{k_i})
		\prod_{i = 0}^\infty \kappa_{i}(d_i).
	\]
	This shows $g_0 \in \overline{G_{\reg}} \cdot \overline{D}$.
\end{proof}

\begin{lemma}\label{AutmGGS:lem:directed elements of G}
	Let $G$ be a non-constant \mGGS-group and let $g \in G$ be an element directed along $\overline 0$. Then $g \in B$.
\end{lemma}

\begin{proof}
	Consider that, since $G \leq \lab(\Sigma)$, there are $(x_1, \dots, x_{p-1}) \in \F_p^{p-1}$ such that
	\[
		\psi_1(g) = (g|_0, a^{x_1}, \dots, a^{x_{p-1}}).
	\]
	Since directed elements stabilise the first layer, there exist $B$-coordinates $\mathbf{n}_{0}, \dots, \mathbf{n}_{p-1}$ and $y_0, \dots, y_{p-1} \in G'$ for $g$. The equation above shows that $\mathbf{n}_1 = \dots = \mathbf{n}_{p-1} = \mathbf 0$ and $y_1 = \dots = y_{p-1} = \id$. Thus the forced $A$-coordinate at $0$ fulfils
	\[
		s_0 = \sum_{i=0}^{p-1} \mathbf n_i \cdot \mathbf c_{p-i} = 0,
	\]
	hence $g|_0 = a^{s_0}b^{\mathbf{n}_{0}} = b^{\mathbf{n}_{0}} \in B$, and in consequence $g = b^{\mathbf{n}_0}\psi_1^{-1}(y_0, \id, \dots, \id)$. Since the set of elements directed along $\overline{0}$ forms a subgroup, the element $\psi_1^{-1}(y_0, \id, \dots, \id)$, and consequently also $y_0$ is directed along $\overline{0}$. We can argue as above for $y_0$, but since $y_0 \in G'$, by \cref{AutmGGS:thm:facts}~(ii), the sum $\sum_{i = 0}^{p-1} \mathbf{n}_i = \mathbf{0}$. Thus $\mathbf{n}_{0} = \mathbf{0}$, and, chasing down the spine, we find $y_0 = \id$. Thus $g \in B$.
\end{proof}

\begin{lemma}\label{AutmGGS:lem:first two terms of overline g1}
	Let $G$ be a non-constant \mGGS-group, and let $h \in \overline{G_{\reg}}$. Then
	\[
		a^h \equiv_{G'} a.
	\]
\end{lemma}

\begin{proof}
	Let $(g_i)_{i\in \N}$ be a sequence of elements $g_i \in G_{\reg}$ for $i \in \N$ such that
	\[
		h = \prod_{i=0}^\infty \kappa_i(g_i) = g_0 \prod_{i=1}^\infty \kappa_i(g_i) = g_0 \kappa_1\left(\prod_{i=0}^\infty \kappa_i(g_{i+1})\right).
	\]
	By \cref{AutmGGS:lem:norm a}, the element $\kappa_1\left(\prod_{i=0}^\infty \kappa_i(g_{i+1})\right)$ centralises $a$. Thus it is sufficient to consider $h = g_0$. The statement now follows from \cref{AutmGGS:lem:facts implications for Greg}~(iii).
\end{proof}

\begin{lemma}\label{AutmGGS:lem:segmenting}
	Let $G$ be a non-constant \mGGS-group. Then 
	\[
		\Nor(G) = \Nor_{\overline{G_{\reg}}}(G) \rtimes \Nor_{\overline{D}}(G)
	\]
\end{lemma}

\begin{proof}
	Assume that $\Nor(G)$ is equal to the product set $\Nor_{\overline{G_{\reg}}}(G) \cdot \Nor_{\overline{D}}(G)$. By \cref{AutmGGS:prop:norm in g1}, $\Nor_{\overline{G_{\reg}}}(G)$ is equal to $G \rtimes \kappa_\infty(A)$ or to $G_{\reg}$. Both groups are normalised by $\Nor_{\overline{D}}(G)$, the first one since $\overline{D}$ normalises $\kappa_\infty(A)$, and the second one since for every $d_0 \kappa_1(d_1)$ with $d_0 \in D$ and $d_1 \in \overline{D}$,
	\[
		\underline{c}^{d_0\kappa_1(d_1)} = \psi_1^{-1}([b,a]^{d_1}, \id, \dots, \id)^{d_0} \in \psi_1^{-1}(G' \dottimes{p} G') \leq G_{\reg}.
	\]
	Thus the product set is in fact a semidirect product. It remains to show the equality $\Nor(G) = \Nor_{\overline{G_{\reg}}}(G) \cdot \Nor_{\overline{D}}(G)$.
	
	By \cref{AutmGGS:lem:norm in ginfty kappa closure}, we may write $g \in \Nor(G)$ as a product $g = h' \cdot d$ with $h' \in \overline{G_{\reg}}$ and $d \in \overline{D}$. Clearly $\overline{D}$ normalises $A$. Thus it is enough to prove: For all $\mathbf{n} \in \F_p^r$ such that $(b^{\mathbf n})^d \not\in G$, then ${h'd} \not\in \Nor(G)$ for all $h' \in \overline{G_{\reg}}$. We to prove that $(h'd)^{-1} \notin \Nor(G)$. Since $\overline{D}$ is a group, we may replace $d$ by its inverse. Write $h$ for $h'^{-1}$. Notice that
	\[
		h = \kappa_1(h_1) h_0
	\]
	for some $h_1 \in \overline{G_{\reg}}$ and $h_0 \in G_{\reg}$. Since $G_{\reg}$ normalises $G$, we may assume that $h_0 = \id$, and thus $h|^\varnothing = \id$. Let $d = \prod_{i = 0}^\infty \kappa_i(d_i)$ for a sequence $(d_i)_{i \in \N}$ of elements $d_i \in D$ such that, for all $i \in \N$, we have $a^{d_i} = a^{j_i}$ for some $j_i \in \Z$. Then, for all $\mathbf{n} \in \F_p^r$,
	\begin{align*}
		\psi_1((b^{\mathbf n})^d) &= ((b^{\mathbf n})^{d|_0}, (a^{\mathbf n \cdot \mathbf{e}_{{1}}})^{d|_0}, \dots, (a^{\mathbf n \cdot \mathbf{e}_{{p-1}}})^{d|_0})^{d_0}\\
		&= ((b^{\mathbf n})^{d|_0}, a^{j_1 \cdot \mathbf n \cdot \mathbf{e}_{1.d_0}}, \dots, a^{j_1 \cdot \mathbf n \cdot \mathbf{e}_{(p-1).d_0}}).
	\end{align*}
	Write $x_{1,i} = j_1 \cdot \mathbf n \cdot \mathbf{e}_{i.d_0}$ for all $i \in \{1, \dots, p-1\}$. Since $\overline{D}$ is self-similar, we see that $(b^{\mathbf n})^d$ is directed along $\overline{0}$, and we write $\mathbf{x}_{k} = (x_{k, 1}, \dots x_{k, p-1})$ for the $A$-exponents of the sections at $i \in \{1, \dots, p-1\}$ of $(b^{\mathbf n})^d|_{0^{k-1}}$.
	
	Now, using \cref{AutmGGS:lem:directed elements of G}, we see that if $(b^{\mathbf n})^d \in G$, then actually $(b^{\mathbf n})^d \in B$. Assume that the $\overline{0}$-directed element $(b^{\mathbf n})^d$ is not a member of $B$. Then there are two possibilities:
	\begin{enumerate}
		\item there exists some $k \in \N$ such that $\mathbf{x}_{k}$ is not contained in the row space of $E$, or,
		\item if for all $k \in \N$ the vector $\mathbf{x}$ is contained in the row space of $E$, but there exists some $k \in \N$ such that $\mathbf{x}_{k} \neq \mathbf{x}_{k+1}$.
	\end{enumerate}
	In both cases, we may assume that $k = 0$, since $\overline{G_{\reg}}$ and $\overline{D}$ are self-similar.

	Given $(b^{\mathbf{n}})^d$, we compute the conjugate by $dh$,
	\begin{align*}
		\psi_1((b^{\mathbf{n}})^{dh}) &= (b^{\mathbf{n}})^d|_0, a^{x_1}, \dots, a^{x_{p-1}})^{h}\\
		&= (((b^{\mathbf{n}})^{d}|_0)^{h|_0}, (a^{x_1})^{h|_1}, \dots, (a^{x_{p-1}})^{h|_{p-1}}).
	\end{align*}
	Since $\overline{G_{\reg}}^{-1}$ is self-similar, we may apply \cref{AutmGGS:lem:first two terms of overline g1}, and we find
	\begin{equation*}\label{eq:dagger}
		\psi_1((b^{\mathbf n})^{dh}) \equiv_{G' \times \dots \times G'} (((b^{\mathbf n})^{d}|_0)^{h|_0} \bmod{G'}, a^{x_1}, \dots, a^{x_{p-1}}).\tag{$\ast$}
	\end{equation*} 
	Assume that we are in case (i), i.e.\ that $\mathbf{x}_{0}$ is not contained in the row space of $E$. Then, by (\ref{eq:dagger}) and \cref{AutmGGS:lem:global equations}, also $(b^{\mathbf n})^{dh} \notin G$.
	
	Assume that we are in case (ii), i.e.\ that $\mathbf{x}_{0} \neq \mathbf{x_{1}}$, but both represent the forced $a$-exponents of an element $b^{\mathbf{m}_{0}}$ and $b^{\mathbf{m}_{1}}$, respectively. Thus by (\ref{eq:dagger}) and \cref{AutmGGS:lem:global equations},
	\begin{align*}
		(b^{\mathbf n})^{dh} &\equiv_{\psi_1^{-1}(G' \times \dots \times G')}  (b^{\mathbf{m}_{0}}) \quad\text{and}\\
		(b^{\mathbf n})^{dh}|_{0} &\equiv_{\psi_1^{-1}(G' \times \dots \times G')}  (b^{\mathbf{m}_{1}})^{h|^{0}}.
	\end{align*}
	Thus
	\[
		(b^{\mathbf{m}_{1}})^{h|^{0}} \equiv (b^{\mathbf n})^{dh}|_{0} \equiv b^{\mathbf{m}_{0}}|_0 \equiv b^{\mathbf{m}_{0}} \bmod{\psi_1^{-1}(G' \times \dots \times G')}.
	\]
	Since $\psi_1^{-1}(G' \times \dots \times G') \cap G = \Stab_G(1)'$, this implies $(b^{\mathbf{m}_{1}})^{a^k} \equiv_{\Stab_G(1)'} b^{\mathbf{m}_{0}}$ for some $k \in \Z$, hence $b^{\mathbf{m}_{1}} = b^{\mathbf{m}_{0}}$ and $\mathbf{m}_{0} = \mathbf{m}_{1}$. But then $\mathbf{x}_{0} = \mathbf{x}_{1}$, a contradiction.
\end{proof}

\section{Elements normalising \(G\) with labels in \(\Delta\)}\label{sec:coprime autom}
Recall that the permutation group $\Delta = \langle \delta \rangle$ is isomorphic to $\F_p^\times$. The rooted automorphism $d= \rooted{\delta}$ acts in two different ways on $G = \Stab_G(1)\rtimes A$. It raises $a$ to a power, i.e.\ it acts my multiplication on the exponent of $a$; and it acts on an element of $g \in \Stab_G(1)$ by permuting the tuple $\psi_1(g)$, i.e.\ by multiplication of the indices of said tuple. Note that the vertex $0$ is fixed by $\delta$.

Recall that $B$ is isomorphic to $\mathbf{E} \leq \F_p^{p-1}$. We now show that $B$ is normalised by every normaliser of $G$ in $\overline{D}$.
\begin{lemma}
	Let $G$ be a non-constant \mGGS-group. Then $\Nor_{\overline{D}}(G) = \Nor_{\overline{D}}(B)$.
\end{lemma}

\begin{proof}
	Since $\overline{D} \leq \Nor(A)$, the inclusion $\Nor_{\overline{D}}(G) \geq \Nor_{\overline{D}}(B)$ is obvious. We now prove the other inclusion. Let $g \in \Nor_{\overline{D}}(G)$. By \cite[Lemma 3.4]{Pet19}, there exists an integer $k \in \N$ such that $g|_{0^k}$ normalises $B$. Thus it is enougth to prove that if $g|_0$ normalises $B$, also $g$ normalises $B$.
	
	Assume $g|_0 \in \Nor(B)$, and let $\mathbf{m}$ and $\mathbf{n} \in \F_p^r$ be such that $(b^{\mathbf{n}})^{g|_0} = b^{\mathbf{m}}$. We may write $g = \kappa_1(g|_0) g|^\varnothing$, where $g|^\varnothing$ normalises $a|^\varnothing$. Hence there exists $j \in \Z$ such that $a^{g|_0} = a^j$. We calculate
	\begin{align*}
		(b^{\mathbf{n}})^{-1}(b^{\mathbf{m}})^{g} &= (b^{\mathbf{n}})^{-1}\psi_1^{-1}((b^{\mathbf{m}})^{g|_0}, a^{j \cdot \mathbf{m} \cdot \mathbf{e}_{{1}}}, \dots, a^{j \cdot \mathbf{m} \cdot \mathbf{e}_{{p_1}}})^{g|^\varnothing}\\
		&= (\id, a^{-\mathbf{n} \cdot \mathbf{e}_{{1}} + j \cdot \mathbf{m} \cdot \mathbf{e}_{{1.g|^\varnothing}}}, \dots, a^{-\mathbf{n} \cdot \mathbf{e}_{{p-1}} + j \cdot \mathbf{m} \cdot \mathbf{e}_{{(p-1).g|^\varnothing}}}).
	\end{align*}
	We see that the commutator coordinates of $(b^{\mathbf{n}})^{-1}(b^{\mathbf{m}})^{g}$ are trivial, the $B$-coordinates are all zero,
	and hence the forced $A$-coordinates are also $s_k = \sum_{i = 0}^{p-1} \mathbf n_i \cdot \mathbf c_{k-i} = 0$. Thus $b^{\mathbf n} = b^{\mathbf m}$.
\end{proof}
Thus, we may restrict our attention to the group $B$. It is fruitful to consider $B$ as a subgroup of the directed subgroup $b^{\F_p^{p-1}}$ of the \mGGS-group associated to the full space $\F_p^{p-1}$ (with standard basis), which is, by the previous lemma, also invariant under $\Nor_{\overline{D}}(G)$.
Write $\mu \colon D \to \F_p^\times$ for the isomorphism induced by
\(
	\alpha^{\delta} = \alpha^{\delta^{\mu}},
\)
where the second operation is taking the power, and define a map $P_\bullet \colon D \to \GL_p(p-1)$ such that $P_d$ is the permutation matrix corresponding to the permutation $d|^\varnothing = \delta$. Let $g = \prod_{i = 0}^\infty \kappa_i(d_i) \in \overline{D}$ for a sequence $(d_i)_{i \in \N}$ of elements in $D$. Then, for all $j \in \{1, \dots, p-1\} = \F_p^\times$,
\begin{align*}
	\psi_1((b^{\mathbf{s}_j})^g) &= ((b^{\mathbf{s}_j})^{g|_0}, \id, \dots, \id, a^{d_1}, \id, \dots, \id)^{d_0}\\
	&= \psi_1((b^{\mathbf{s}_j})^{g|_0}, \id, \dots, \id, (d_1)^\mu a, \id, \dots, \id),
\end{align*}
where the non-trivial entries (in the second line) are the positions $0$ and $(d_1)^\mu j$. If $g$ normalises $G$, it must normalise $B$, and the conjugate $(b^{\mathbf{s}_j})^g$ is determined by the exponents of the sections at the positions in $\{1, \dots, p-1\}$. Thus
\[
	(b^{\mathbf{s}_j})^g = b^{(d_1)^\mu \mathbf{s}_{(d_0)^{\mu} j}}.
\]
In other words, the action of $g \in \Nor_{\overline{D}}(G)$ induces, via the isomorphism $b^{\bullet}$, the linear map
\begin{equation*}\label{eq:linear map}
	(d_1)^{\mu} P_{d_0}\tag{\textdagger}
\end{equation*}
on $\F_p^{p-1}$. Returning to the directed group $B$, we see that every $g \in \Nor_{\overline{D}}(G)$ must be such that $P_{d_0}$ leaves $\mathbf{E}$ invariant. Hence we define
\[
	U \defeq \{ u \in D \mid \mathbf{E}P_u = \mathbf{E} \} = \stab_D(\mathbf E).
\]
Furthermore, we define the subgroup
\begin{align*}
	V \defeq& \{ v \in U \mid \text{ for all } \mathbf{e} \in \mathbf{E} \text{ there exists } \lambda \in \F_p^\times \text{ such that } \mathbf{e}P_v = \lambda \mathbf{e} \} \\
	=& \{ v \in U \mid \mathbf{E} \subseteq \operatorname{Eig}_\lambda(P_v) \text{ for some }\lambda \in \F_p^\times \}.
\end{align*}
Since $V \leq D$ is cyclic, the element $\lambda \in \F_p^\times$ generating a maximal subgrop is uniquely determined. Finally, define the subgroup
\[
	W \defeq \langle \lambda \rangle.
\]
\begin{proposition}\label{AutmGGS:prop:norm in D}
	Let $G$ be a non-constant \mGGS-group. Let $U$ and $W$ be defined as above. Then
	\begin{align*}
		\Nor_{\overline{D}}(G)
		\cong U \times W.
	\end{align*}
\end{proposition}

\begin{proof}
	By (\ref{eq:linear map}), the action of a given element $g = \prod_{i = 0}^\infty \kappa_{i}(d_i)$ is determined by $d_0$ and $d_1$. Furthermore, since $\mathbf{E}$ must be invariant under $P_{d_0}$, we see that necessarily $d_0 \in U$. Since $\Nor(G)$, by \cite[Lemma~3.3]{Pet19}, and $\overline{D}$ are self-similar, we find, for all $k \in \N$,
	\[
		g|_{0^k} = \prod_{i = 0}^\infty \kappa_{i}(d_{i+k}) \in \Nor_{\overline{D}}(G).
	\]
	Since, for all $\mathbf{n} \in \F_p^r$ and $g \in \Nor_{\overline{D}}(G)$,
	\[
		(b^{\mathbf{n}})^g = (b^{\mathbf{n}})^g|_0 = (b^{\mathbf{n}}|_0)^{g|_0} = (b^{\mathbf{n}})^{g|_0},
	\]
	we see that the action induced on $\F_p^{p-1}$ by all elements $g|_{0^k}$, for $k \in \N$, are equal, i.e.\ that the following equalities of matrices hold,
	\[
		(d_{k+1})^{\mu} P_{d_k} = (d_1)^{\mu} P_{d_0}, \quad\text{ hence }\quad \mathrm{I}_{p-1} = (d_{k+1}d_{k+2}^{-1})^{\mu} P_{d_k d_{k+1}^{-1}},
	\]
	where $\mathrm{I}_{p-1}$ is the identity matrix. Recall that, for all $k \in \N$, the matrix $P_{d_k d_{k+1}^{-1}}$ acts either does not act as a scalar on $\mathbf{E}$, hence there is no $d_{k+1}$ fulfilling the equation above, or it acts as some scalar $\lambda^i$, for some $i \in \Z$. Thus, every difference $d_kd_{k+1}^{-1}$ must be an element of $W$, otherwise, $g$ cannot be normalising $G$.
	
	On the other hand, for $d_0 \in U$ and $d_1 \in d_0^{-1}W$ there is a unique sequence $(d_i)_{i \in \N}$ that defines an element of $\Nor_{\overline{D}}(G)$, since
	\begin{align*}
		d_{k+2} &= d_{k+1} (d_kd_{k+1}^{-1})^{r}\\
		&= d_{k+1} (d_k ((d_{k-1}d_k^{-1})^r)^{-1} d_k^{-1})^r\\
		&= d_{k+1} ((d_{k-1}d_k^{-1})^{-1})^{r^{2}}\\
		&= d_{k+1} ((d_0 d_1^{-1})^{(-1)^k})^{r^{k+1}}.
	\end{align*}
	Thus $\Nor_{\overline{D}}(G) \cong U \times W$.
\end{proof}
In particular, if $r = 1$, every linear map leaving $\mathbf{E}$ invariant is a scalar multiplication, i.e.\ the subgroups $U$ and $V$ coincide.
Clearly, if $W$ is the trivial group, the only elements of $\Nor_M(G)$ are defined by the constant sequences. More generally, the sequence $(d_i)_{i\in \N}$ defining the normalising element with given $d_0$ and $d_1$ is periodic with periodicity prescribed by the order of $\lambda$.

Now all ingredients are ready for the proof of our main theorem.

\begin{proof}[Proof of \cref{AutmGGS:thm:main}]
	By \cite[Theorem~1]{GW03} and \cite[Proposition~3.7]{KT18}, the automorphism group of $G$ coincides with the normaliser of $G$ in $\Aut(X^*)$. By \cref{AutmGGS:lem:segmenting}, this normaliser is the semidirect product
	\[
		\Nor_{\overline{G_{\reg}}}(G) \rtimes \Nor_{\overline{D}}(G).
	\]
	These two groups were computed in \cref{AutmGGS:prop:norm in g1} and \cref{AutmGGS:prop:norm in D}.
\end{proof}

\section{Examples}\label{sec:examples automggs}

To illustrate the definitions of $U, V$ and $W$ we compute some explicit examples.

\begin{example}
	Let $G$ be the \GGS-group acting on the $5$-adic tree with $\mathbf{E}$ generated by $(1,2,2,1)$. Clearly $G$ is symmetric. For every symmetric \GGS-group, the space $\mathbf{E}$ is by definition invariant under the permutation induced by $-1 \in \F_p^{\times}$. In fact, it always acts trivially, hence $-1 \in W$. In our case, this is the only non-trivial permutation leaving $\mathbf{E}$ invariant, since
	\[
		(1,2,2,1) P_{x \mapsto 2x} = (2,1,1,2) = (1,2,2,1) P_{x \mapsto 3x},
	\]
	while $(2,1,1,2) \notin \mathbf{E}$. Thus
	\[
		\Aut(G) = (G \rtimes \langle \underline{c} \rangle) \rtimes \langle \textstyle \prod_{i = 0}^\infty \kappa_{i}(x \mapsto -x) \rangle \cong (G \rtimes \mathrm{C}_5) \rtimes \mathrm{C}_2.
	\]
	The group $G$ and the group defined by $(1,4,4,1)$ are the \mGGS-groups with the smallest possible outer automorphism group.
\end{example}

\begin{example}
	Let $G$ be the (regular) \GGS-group acting on the $p$-adic tree with $E$ generated by $\mathbf{b} = (1,2,\dots,p-1)$. Let $(\lambda_1, \dots, \lambda_{p-1})$ be the image of $\mathbf{b}$ under $P_d$, for $d \in D$. Since
	\[
		\lambda_i = \mathbf{b}_{d^{-1}i} = (d^{-1})^\mu i = (d^{-1})^\mu \mathbf{b}_{i},
	\]
	we see that $\mathbf{b}P_d = (d^{-1})^\mu \mathbf{b}$. Thus $W = V = U = D$, and the automorphism group is `maximal',
	\[
		\Aut(G) = (G \rtimes \kappa_{\infty}(A)) \rtimes (\F_p^\times)^2.
	\]
\end{example}

\begin{example}
	The distinction between the subgroups $U, V$ and $W$ is not superficial. Consider the vector $\mathbf{b}_1 = (1, 2, 11, 3, 12, 10, 10, 12, 3, 11, 2, 1) \in \F_{13}^{12}$. An easy calculation shows that
	\[
		\mathbf{b}_1 P_{x \mapsto 5x} = (12, 11, 2, 10, 1, 3, 3, 1, 10, 2, 11, 12) = -\mathbf{b}_1,
	\]
	while $\mathbf{b}_1 P_{x \mapsto 3x}$ is not a multiple of $\mathbf{b}_1$. Set $\mathbf{b}_2 = \mathbf{b}_1 P_{x \mapsto 3x}$ and $\mathbf{b}_3 = \mathbf{b}_1 P_{x \mapsto 9x}$ and let $\mathbf{E}$ be the space spanned by $\mathbf{b}_1, \mathbf{b}_2$ and $\mathbf{b}_3$. Since $\F_{13}^\times$ is generated by $3$ and $5$, the space $\mathbf{E}$ is invariant under all permutations induced by index-multiplication, i.e.\ $U = D$. But only the multiplication by the multiples of $5$ act by scalar multiplication of $\mathbf{E}$, hence $V = \langle x \mapsto 5x \rangle$. The corresponding scalars are $1$ and $12$, hence $W$ is of order $2$.
\end{example}

\begin{example}
	Let $G_p$ be a Gupta--Sidki $p$-group, i.e.\ the \GGS-group with $\mathbf{E}$ spanned by $\mathbf{b} = (1, -1, 0, \dots, 0) \in \F_p^{p-1}$. All Gupta--Sidki $p$-groups are regular. Let $n \in N$, and consider $\mathbf{b}P_n$. Since the projection to the last $p-3$ coordinates of $\mathbf{E}$ is trivial, the index $1$ must be mapped to $1$ or $2$ under $n$, and the same holds for $2$. This is only possible if $n = 1$ or $n = 2$ and $2\cdot 2 \equiv_p 1$, hence in case $p = 3$. Otherwise $U$ is trivial. If $p = 3$, the group $W$ is equal to $U$, since the non-trivial permutation induced by the index multiplication by $2$ is equal to pointwise multiplication by $2$. This recovers the result of \cite{Sid87}, where the automorphism group of $G_3$ was first computed. Interestingly, this example is the `odd one out', having automorphisms of order~$2$.

	Concludingly, we found
	\[
		\Aut(G_p) = \begin{cases}
			(G_p \rtimes \kappa_{\infty}(A)) \rtimes \mathrm{C}_2^2 &\text{ if }p = 3,\\
			G_p \rtimes \kappa_{\infty}(A) &\text{ otherwise. }
		\end{cases}
	\]
\end{example}

\begin{example}
	Let $G_{\F_p^{p-1}}$ be the \mGGS-group defined by the full space $\F_p^{p-1}$. This group is regular, and every permutation $P_n$ leaves $\F_p^{p-1}$ invariant. On the other hand, no non-trival permutation acts on the full space as a multiplication. Thus 
	\begin{align*}
		\Aut(G_{\F_p^{p-1}}) &= (G_{\F_p^{p-1}} \rtimes \kappa_{\infty}(A)) \rtimes \{ \textstyle \prod_{i = 0}^\infty \kappa_{i}(d') \mid d' \in D\}\\
		&\cong (G_{\F_p^{p-1}} \rtimes \prod_{\omega} \mathrm{C}_p) \rtimes \mathrm{C}_{p-1}.
	\end{align*}
\end{example}


\begin{thebibliography}{10}

\bibitem{AKT16}
T.~Alexoudas, B.~Klopsch, and A.~Thillaisundaram.
\newblock Maximal subgroups of multi-edge spinal groups.
\newblock {\em Groups, Geometry, and Dynamics}, 10(2):619--648, 2016.

\bibitem{BS05}
L.~Bartholdi and S.~N. Sidki.
\newblock The automorphism tower of groups acting on rooted trees.
\newblock {\em Transactions of the American Mathematical Society},
  358(1):329--358, Mar. 2005.

\bibitem{FZ13}
G.~{Fern{\'a}ndez-Alcober} and A.~{Zugadi-Reizabal}.
\newblock {{GGS-groups}}: {{Order}} of congruence quotients and {{Hausdorff}}
  dimension.
\newblock {\em Transactions of the American Mathematical Society},
  366(4):1993--2017, Oct. 2013.

\bibitem{FGU17}
G.~A. {Fern{\'a}ndez-Alcober}, A.~Garrido, and J.~{Uria-Albizuri}.
\newblock On the congruence subgroup property for {{GGS-groups}}.
\newblock {\em Proceedings of the American Mathematical Society},
  145(8):3311--3322, Jan. 2017.

\bibitem{GU19}
A.~Garrido and J.~{Uria-Albizuri}.
\newblock Pro-{{{\emph{C}}}} congruence properties for groups of rooted tree
  automorphisms.
\newblock {\em Archiv der Mathematik}, 112(2):123--137, Feb. 2019.

\bibitem{Gri80}
R.~I. Grigorchuk.
\newblock {On Burnside's problem on periodic groups}.
\newblock {\em Funktsional. Anal. i Prilozhen.}, 14(1):53--54, 1980.

\bibitem{GS04}
R.~I. Grigorchuk and S.~N. Sidki.
\newblock The group of automorphisms of a 3-generated 2-group of intermediate
  growth.
\newblock {\em International Journal of Algebra and Computation},
  14(05n06):667--676, Oct. 2004.

\bibitem{GW03}
R.~I. Grigorchuk and J.~S. Wilson.
\newblock The uniqueness of the actions of certain branch groups on rooted
  trees.
\newblock {\em Geometriae Dedicata}, 100(1):103--116, 2003.

\bibitem{GS83}
N.~Gupta and S.~Sidki.
\newblock On the {{Burnside}} problem for periodic groups.
\newblock {\em Mathematische Zeitschrift}, 182(3):385--388, Sept. 1983.

\bibitem{KT18}
B.~Klopsch and A.~Thillaisundaram.
\newblock Maximal {{Subgroups}} and {{Irreducible Representations}} of
  {{Generalized Multi-Edge Spinal Groups}}.
\newblock {\em Proceedings of the Edinburgh Mathematical Society},
  61(3):673--703, Aug. 2018.

\bibitem{LN02}
Y.~Lavreniuk and V.~Nekrashevych.
\newblock Rigidity of branch groups acting on rooted trees.
\newblock {\em Geometriae Dedicata}, 89(1):155--175, 2002.

\bibitem{Pet19}
J.~M. Petschick.
\newblock On conjugacy of {{GGS-groups}}.
\newblock {\em Journal of Group Theory}, 22(3):347--358, May 2019.

\bibitem{Rov02}
C.~E. R{\"o}ver.
\newblock Abstract commensurators of groups acting on rooted trees.
\newblock {\em Geometriae Dedicata}, 94(1):45--61, 2002.

\bibitem{Sid87}
S.~Sidki.
\newblock On a 2-generated infinite 3-group: {{Subgroups}} and automorphisms.
\newblock {\em Journal of Algebra}, 110(1):24--55, Oct. 1987.

\end{thebibliography}

\end{document}